\documentclass[11pt]{article}

\newtheorem{Def}{Definition}[section]
\newtheorem{Th}{Theorem}[section]
\newtheorem{Prop}{Proposition}[section]
\newtheorem{Not}{Notation}[section]
\newtheorem{Lemma}{Lemma}[section]
\newtheorem{Rem}{Remark}[section]
\newtheorem{Cor}{Corollary}[section]
\newtheorem{Conj}{Conjecture}[section]

\def\s{\section}
\def\ss{\subsection}

\def\d{\begin{Def}}
\def\t{\begin{Th}}
\def\p{\begin{Prop}}
\def\n{\begin{Not}}
\def\la{\begin{Lemma}}
\def\r{\begin{Rem}}
\def\c{\begin{Cor}}
\def\cj{\begin{Conj}}
\def\ee{\begin{equation}}
\def\aa{\begin{eqnarray}}
\def\y{\begin{eqnarray*}}
\def\bd{\begin{description}}

\def\ed{\end{Def}}
\def\et{\end{Th}}
\def\ep{\end{Prop}}
\def\en{\end{Not}}
\def\el{\end{Lemma}}
\def\er{\end{Rem}}
\def\ec{\end{Cor}}
\def\ecj{\end{Conj}}
\def\eee{\end{equation}}
\def\eaa{\end{eqnarray}}
\def\ey{\end{eqnarray*}}
\def\ebd{\end{description}}

\def\nn{\nonumber}
\def\bp{{\bf Proof.}\hspace{2mm}}
\def\qe{\hfill{\rm Q.E.D.}}

\def\ox{\mbox{}}
\def\lb{\label}
\def\bs{\setminus}

\def\R{{\bf R}}
\def\C{{\bf C}}
\def\Z{{\bf Z}}
\def\N{{\bf N}}
\def\Ua{{\bf U}}

\def\Ha{{\cal H}}
\def\SH{{\cal SH}}

\def\Pa{{\cal P}}
\def\J{{\cal J}}

\def\I{{\cal I}}
\def\E{{\cal E}}

\def\O{{\cal O}}

\def\gl{{\rm gl}}
\def\GL{{\rm GL}}
\def\Sp{{\rm Sp}}

\def\vp{\epsilon}
\def\im{{\rm im}}
\def\ind{{\rm ind}}

\title{Multiplicity of closed characteristics on symmetric
convex hypersurfaces in $\R^{2n}$}
\author{Chun-gen Liu\thanks{Partially supported by NNSF(10071040), Exc. Ph.D.
Funds of ME of China, PMC Key Lab of EM of China, and the Hong Kong Qiu Shi
Sci. Tech. Foundation.}\\
Department of Mathematics, Nankai University\\
Tianjin 300071, The People's Republic of China\\
\\   and  \\
Yiming Long,\thanks{Partially supported by the 973 Program of STM,
NNSF, MCME, RFDP, PMC Key Lab of EM of China, S. S. Chern Foundation,
Hong Kong Qiu Shi Sci. Tech. Foundation, and CEC of Tianjin. }
\thanks{Associate member of ICTP.}
\qquad Chaofeng Zhu\thanks{Partially supported by the Hong Kong Qiu Shi
Sci. Tech. Foundation.}\\
Nankai Institute of Mathematics, Nankai University\\
Tianjin 300071, The People's Republic of China}
\date{}

\begin{document}

\maketitle

\begin{abstract}
{\it Let $\Sigma$ be a compact $C^2$ hypersurface in $\R^{2n}$
bounding a convex set with non-empty interior.
In this paper it is proved that there always exist at least
$n$ geometrically distinct closed characteristics on $\Sigma$ if
$\Sigma$ is symmetric with respect to the origin.}
\end{abstract}

\s{Introduction and main results.}\lb{99mls1}

Our aim in this paper is to study the multiplicity of closed characteristics
on any $C^2$-convex compact smooth hypersurface in $\R^{2n}$ which is
symmetric with respect to the origin. Let $\Sigma$ be a $C^2$-compact
hypersurface in $\R^{2n}$ bounding a convex compact set $C$ with non-empty
interior, possess a non-vanishing Gaussian curvature, and $0\in C$.
We denote the set of all such hypersurfaces in $\R^{2n}$ by
$\Ha(2n)$, and the set $\{\Sigma\in\Ha(2n)\mid\Sigma=-\Sigma\}$ by
$\SH(2n)$, where $-\Sigma=\{x\in\R^{2n}\mid -x\in\Sigma\}$.
For $x\in\Sigma$, let $N_{\Sigma}(x)$ be the outward unit normal vector at
$x$ on $\Sigma$. We consider the given energy problem of finding $\tau>0$
and an absolutely continuous curve $x\colon[0,\tau]\to\R^{2n}$ such that
\ee \cases{\dot{x}(t)&$= JN_{\Sigma}(x(t)), \quad x(t)\in\Sigma,
\qquad\forall t\in\R,$\cr
x(\tau) &$= x(0),$ \cr} \lb{99che101}\eee
where $J=\pmatrix{0&-I_n\cr I_n&0\cr}$, $I_n$ is the identity matrix on
$\R^n$. When there is no confusion we shall omit the subindex of
the identity matrices. A solution $(\tau,x)$ of the problem (\ref{99che101})
is called a {\bf closed characteristic} on $\Sigma$. Two closed
characteristics $(\tau,x)$ and $(\sigma,y)$ are {\bf geometrically distinct},
if $x(\R)\not= y(\R)$. We denote by $\J(\Sigma)$ and $\tilde\J(\Sigma)$ the
set of all closed characteristics $(\tau,x)$ on $\Sigma$ with $\tau$ being
the minimal period of $x$ and the set of all geometrically distinct ones
respectively. For $(\tau,x)\in\J(\Sigma)$, we denote by $[(\tau,x)]$ the set
of all elements in $\J(\Sigma)$ which is geometrically the same as $(\tau,x)$.
$^{\#}A$ denotes the total number of elements in a set $A$.

The study on closed characteristics in the global sense started in 1978,
when the existence of at least one closed characteristic on any
$\Sigma\in\Ha(2n)$ was first established by P. Rabinowitz in \cite{Ra1}
(for star-shaped hypersurfaces) and A. Weinstein in \cite{We2} independently.
In \cite{EH} of I. Ekeland and L. Lassoued, \cite{EL} of I. Ekeland and H.
Hofer in 1987, and \cite{Sz} of A. Szulkin in 1988,
$^{\#}\tilde\J(\Sigma)\ge 2$ was proved for any $\Sigma\in\Ha(2n)$ when
$n\ge 2$. In the recent paper \cite{LZh2}, Y. Long and C. Zhu proved
$^{\#}\tilde\J(\Sigma)\ge[\frac{n}{2}]+1$ for any $\Sigma\in\Ha(2n)$,
where $[a]$ denotes the greatest integer which is not greater than $a$.
On the other hand, $^{\#}\tilde\J(\Sigma)\ge n$ was proved by I. Ekeland
and J. M. Lasry in \cite{ELr} of 1980 for $\Sigma\in\Ha(2n)$ which is
$\sqrt{2}$-pinched, by M. Girardi in \cite{G} of 1984 for $\Sigma\in\SH(2n)$
which is $\sqrt{3}$-pinched. Other related significant progresses can be
found in \cite{We1} and \cite{Mo} for local results, and in \cite{BLMR},
\cite{Ek3}, and \cite{HWZ} for global results.

A typical example of $\Sigma\in\Ha(2n)$ is the ellipsoid $\E_n(r)$ defined
as follows. Let $r=(r_1,\ldots,r_n)$ with $r_k>0$ for $1\leq k\leq n$.
Define
\ee \E_n(r)=\{x=(x_1,\ldots,x_n)\in\R^{2n}\,|\,
      \frac{1}{2}\sum_{k=1}^n{{|x_k|^2}\over {r_k^2}}=1\}. \lb{99che103}\eee
If $r_j/r_k$ is irrational whenever $j\not=k$, this $\E_n(r)$ is called a
{\bf weakly non-resonant ellipsoid}. In this case, it is well known that
there are precisely $n$ geometrically distinct closed characteristics on
$\E_n(r)$ (cf. \S{I.7} of \cite{Ek3}).

  A long standing conjecture mentioned by I. Ekeland on Page 235 of
\cite{Ek3} is
$$ ^{\#}\tilde\J(\Sigma)\ge n, \qquad \forall \Sigma\in \Ha(2n). $$
Our following main result in this paper gives a positive answer to
this conjecture for $\Sigma\in\SH(2n)$.

\t \lb{99mlt101} For every $\Sigma\in\SH(2n)$, we have
\ee   ^{\#}\tilde\J(\Sigma)\ge n.     \lb{99mcl106}\eee
\et

  For reader's convenience, we briefly sketch the main idea of the proof
below. Fix a $\Sigma\in\SH(2n)$ and suppose
$\,^{\#}\tilde{\J}(\Sigma)<+\infty$. Our proof of this theorem is carried
out in the following $7$ steps.

  $1^{\circ}$ Using the method of I. Ekeland (cf. Section V.3 of \cite{Ek3}),
this problem is transformed to a fixed energy problem, and then to a $1$-periodic
solution problem of a Hamiltonian system, whose Hamiltonian function is defined
by $H_{\alpha}(x)=j_{\Sigma}(x)^{\alpha}$ for any $x\in\R^{2n}$ in terms of
the gauge function $j_{\Sigma}(x)$ of $\Sigma$ and some $\alpha\in (1,2)$.
Solutions of this $1$-periodic solution problem correspond to critical points of
the Clarke-Ekeland dual action functional $f_{\alpha}:E_{\alpha}\to\R$ given
by (\ref{99che302})-(\ref{99che303}) below, and $f_{\alpha}$ possesses a strictly
increasing infinite sequence of critical values $\{c_k\}$ with $c_k\to 0$ as
$k\to +\infty$. These $c_k$'s are obtained by the Fadell-Rabinowitz
$S^1$-cohomological index method as in Section V.3 of \cite{Ek3}. Note that
here each $[(\tau,x)]\in\tilde{\J}(\Sigma)$ and all of its iterations
$\{x^m\}$ yield a strictly increasing infinite subsequence of $\{c_k\}$.
Denote the Maslov-type index interval of the $m$-th iteration of $(\tau,x)$ by
$$  \I_m(\tau,x)=[i_{m\tau}(x^m), i_{m\tau}(x^m)+\nu_{m\tau}(x^m)]. $$
Then $\I_m(\tau,x)\cap\I_{m'}(\tau,x)=\emptyset$ if $m\not= m'$, and
$\{i_{m\tau}(x^m)\}$ is a strictly increasing sequence. The union of all
the Maslov-type index intervals of all iterations of all the elements
in $\tilde{\J}(\Sigma)$ covers the integer set $2\N-2+n$, where $\N$ is the
set of positive integers. Thus instead of studying the arrangements of the
subsequences of $\{c_k\}$ of elements in $\tilde{\J}(\Sigma)$ in the whole
sequence $\{c_k\}$, we study how their Maslov-type index interval sequences
cover the set $2\N-2+n$.

  $2^{\circ}$ Let $\O(x)=x(\R)$ for any closed characteristic
$(\tau,x)\in\J(\Sigma)$. Our Lemma 4.2 below shows that every
$(\tau,x)\in\J(\Sigma)$ is either {\bf symmetric}, i.e. $\O(x)=\O(-x)$,
or {\bf asymmetric}, i.e. $\O(x)\cap\O(-x)=\emptyset$. Moreover, any
symmetric $(\tau,x)$ satisfies
\ee  x(t+\frac{\tau}{2}) = -x(t), \qquad \forall t\in\R. \lb{eq11.1}\eee
any asymmetric $(\tau,x)$ satisfies
\ee  (i_{m\tau}(x^m),\nu_{\tau}(x^m))
   = (i_{m\tau}((-x)^m),\nu_{\tau}((-x)^m)), \quad \forall m\in\N.
       \lb{eq11.1a}\eee
Denote the numbers of symmetric and asymmetric elements in
$\tilde{\J}(\Sigma)$ by $q_1$ and $2q_2$. We can write
\ee
\tilde\J(\Sigma)=\{[(\tau_j,x_j)]\,|\,j=1,\ldots,q_1\}
\cup\{[(\tau_k,x_k)],[(\tau_k,-x_k)]|k=q_1+1,\ldots,q_1+q_2\}.
\lb{eq11.2}\eee

  $3^{\circ}$ Applying the common index jump theorem proved by Y. Long and
C. Zhu in \cite{LZh2} (cf. Theorem 4.1 below) to fundamental solutions
of the linearized Hamiltonian systems at
$$ (\tau_1,x_1),\ldots,(\tau_{q_1+q_2},x_{q_1+q_2}),
  (2\tau_{q_1+1},x_{q_1+1}^2),\ldots,(2\tau_{q_1+q_2},x_{q_1+q_2}^2), $$
we obtain an integer $N$ and iteration times $m_1,\ldots,m_{q_1+2q_2}$ such
that the rather precise information on the Maslov-type indices of
iterations of $(\tau_j,x_j)$'s are given in (\ref{99mle442})-(\ref{99mle449}).

  $4^{\circ}$ Combining above $1^{\circ}$ and $3^{\circ}$ together, by the
property of Fadell-Rabinowitz $S^1$-cohomology index (cf. Lemma 3.1 below),
we obtain an injection map $p:\N\to \tilde{\J}(\Sigma)\times \N$ such that
\ee  p(N-s+1) = ([(\tau_{k(s)},x_{k(s)})],m(s)), \qquad
        {\rm for}\;\;s=1,\ldots,n, \lb{eq11.3}\eee
where $m(s)$ is the iteration time of $(\tau_{k(s)},x_{k(s)})$ such that
$$ 2(N-s+1)-2+n \in \I_{m(s)}(\tau_{k(s)},x_{k(s)}). $$

  $5^{\circ}$ Let
\ee S_1=\{s\in\{1,\ldots,n\}\,|,k(s)\le q_1\}, \quad
   S_2=\{1,\ldots,n\}\setminus S_1. \lb{eq11.4}\eee
In Section 4, we prove
\ee  ^{\#}S_1 \le q_1, \quad {\rm and}\quad ^{\#}S_2 \le 2q_2.
    \lb{eq11.5}\eee
Then this implies
$$ ^{\#}\tilde{\J}(\Sigma) = q_1+2q_2 \ge ^{\#}S_1 + ^{\#}S_2 =n. $$

  $6^{\circ}$ To prove the first estimate in (\ref{eq11.5}), using the
property (\ref{eq11.1}) of symmetric orbit $(\tau,x)$ and precise index
information (\ref{99mle442})-(\ref{99mle449}) of iterations of
$(\tau_j,x_j)$'s, we conclude that the Maslov-type index, nullity,
and splitting numbers of $x|_{[0,\tau]}$ are the same as those of
$(x|_{[0,\tau/2]})^2$. Therefore there holds
$$ i_{\tau}(x) + 2S^+(x)-\nu_{\tau}(x) \ge n. $$
By this estimate, we get that the integer $m(s)$ in (\ref{eq11.3}) is
uniquely determined by $k(s)$ there provided $k(s)\le q_1$. Then the
injection map $p$ of (\ref{eq11.3}) induces an injection map from $S_1$
to $\{[(\tau_j,x_j)]\,|\,1\le j\le q_1\}$ and yields the first estimate
in (\ref{eq11.5}).

  $7^{\circ}$ To prove the second estimate in (\ref{eq11.5}), using the
property (\ref{eq11.1a}) of asymmetric orbit $(\tau,x)$ and precise index
information (\ref{99mle442})-(\ref{99mle449}) of iterations of
$(\tau_j,x_j)$'s, we conclude that the integer $m(s)$ in (\ref{eq11.3})
possesses at most two choices determined by the $k(s)$ there, provided
$q_1< k(s)\le q_1+q_2$. Then the injection map $p$ of (\ref{eq11.3})
induces a map from $S_2$ to
$\Gamma\equiv\{[(\tau_j,x_j)]\,|\,q_1< j\le q_1+q_2\}$ such that any
element in $\Gamma$ is the image of at most two numbers in $S_2$. This
yields the second estimate in (\ref{eq11.5}), and completes the proof
of Theorem \ref{99mlt101}.

  This paper is organized as follows. In \S{2}, we briefly review
the Maslov-type index theory. In \S{3}, we discuss main variational
properties of closed characteristics and give details of the above
step $1^{\circ}$.  In \S{4}, we complete the above steps
$2^{\circ}$-$6^{\circ}$ and prove Theorem \ref{99mlt101}.

\s{The Maslov-type index theory and its iteration theory.}
\lb{99mls2}

\ss{Maslov-type index theory.}\lb{99mlss21}

In this subsection we give a brief review on the Maslov-type index theory for
symplectic matrix paths.

For any $n\in\N$ and $\tau>0$, we define as usual
\y \Sp(2n)&=&\left\{M\in\GL(2n,\R)\mid M^TJM=J\right\},\\
\Pa_{\tau}(2n)&=&\left\{\gamma\in C([0,\tau],\Sp(2n))\,|\,\gamma(0)=I\right\}.
\ey

For any $\gamma\in\Pa_{\tau}(2n)$, the Maslov-type index of $\gamma$ is
defined to be a pair of integers by C. Conley, E. Zehnder, and Y. Long in
the works \cite{CZ}, \cite{LZ}, \cite{Lo1}, and \cite{Lo2}, which is
denoted by
$$ (i_{\tau}(\gamma),\nu_{\tau}(\gamma))\in\Z\times\{0,1,\ldots,2n\}, $$
where $\Z$ is the set of all integers. The iteration theory of this
Maslov-type index theory was established in \cite{Lo4} by Y. Long via the
$\omega$-index theory introduced there. Let $\Ua\equiv\{z\in\C\,|\,|z|=1\}$.
For any $\omega\in\Ua$ and $\gamma\in\Pa_{\tau}(2n)$, the
$\omega$-index of $\gamma$ is denoted by
$$ (i_{\tau,\omega}(\gamma),\nu_{\tau,\omega}(\gamma))
      \in\Z\times\{0,1,\ldots,2n\}. $$
Note that there hold $i_{\tau,1}(\gamma)= i_{\tau}(\gamma)$ and
$\nu_{\tau,1}(\gamma)=\nu_{\tau}(\gamma)$. The details of the definitions
and properties of these index and their iteration theories can be found
in the above mentioned works as well as \cite{Lo5} and \cite{Lo6}. We also
refer to \cite{LZh} for another approach.

For $B\in C(\R/(\tau\Z), \gl(2n,\C))$ with $B(t)$ being self-adjoint for
all $t$, we consider the linear Hamiltonian system
\ee \dot{x}=JB(t)x,\qquad x\in\C^{2n}.\lb{99che203}\eee

The following proposition gives the relationship between the Maslov-type
index theory and the Ekeland index theory.

\p \lb{99chp203} (\cite{Br}, Lemma 1.3 of \cite{Lo3}, and
Theorem 3.2 of \cite{LZh}) For $B\in C(\R/\tau\Z, \gl(2n,\R))$ with $B(t)$
being positively definite and symmetric for all $t$, let
$\gamma_B\in\Pa_{\tau}(2n)$ be the fundamental solution of the linear
Hamiltonian system (\ref{99che203}), $i_{\tau}^E(\gamma_B)$ and
$\nu_{\tau}^E(\gamma_B)$ be the Ekeland index and nullitity given
by Definition I.4.3 in \cite{Ek3} with $J$ replaced by that in
(\ref{99che101}). Then we have
\aa i_{\tau}^E(\gamma_B) &=& i_{\tau}(\gamma_B)-n,\lb{99che208}\\
 \nu_{\tau}^E(\gamma_B)  &=& \nu_{\tau}(\gamma_B).\lb{99che209}\eaa
\ep

\c \lb{99chc201} For any positive definite symmetric path
$B\in C(\R/(\tau\Z),\gl(2n,\R))$, we have
\ee \lb{99che210} i_{\tau}(\gamma_B)\ge n. \eee\ec

\ss{Bott-type formulae for splitting numbers.}\lb{99mlss22}

For $n\in\N$, $\tau>0$ and a path $\gamma\in\Pa_{\tau}(2n)$,
we define the {\bf iteration path} $\tilde\gamma$ of $\gamma$ by
\ee \tilde\gamma(t)=\gamma(t-j\tau)\gamma(\tau)^j,\qquad\mbox{for}\;\,
j\tau\le t\le (j+1)\tau \;\,\mbox{and}\;\, j\in\{0\}\cup\N.\lb{752}\eee

For $M\in\Sp(2n)$, let $S^{\pm}_M(\omega)$
be the splitting numbers defined by Theorem 1.3 of \cite{Lo4} for
all $\omega\in\Ua$. Then the following Bott-type formula follows.

\la\lb{99mll201} For any $M\in\Sp(2n)$, $m\in\N$ and $z\in\Ua$, there holds
\ee S^{\pm}_{M^m}(z)=\sum_{\omega^m=z}S^{\pm}_M(\omega).\lb{99ml201}\eee
\el

\bp Since $\Sp(2n)$ is path-connected, we can choose a path
$\gamma\in\Pa_{\tau}(2n)$ with $\gamma(\tau)=M$. Define
$\theta\in[0,2\pi)$ by $z=e^{\sqrt{-1}\theta}$. Let $\tilde\gamma$ be the
iteration path of $\gamma$ defined by (\ref{752}). By the definition of
the splitting numbers and Theorem 1.4 in \cite{Lo4}, we obtain
\y S^{+}_{M^m}(z)
&=&\lim_{\vp\to 0^+}i_{m\tau,e^{\sqrt{-1}(\theta+\vp)}}(\tilde\gamma)
-i_{m\tau,e^{\sqrt{-1}\theta}}(\tilde\gamma)\\
&=&\lim_{\vp\to 0^+}\sum_{k=1}^m
i_{\tau,e^{\sqrt{-1}(\theta+\vp+2k\pi)/m}}(\gamma)
-\sum_{k=1}^m i_{\tau,e^{\sqrt{-1}(\theta+2k\pi)/m}}(\gamma)\\
&=&\sum_{k=1}^m\left(\lim_{\vp\to 0^+}
i_{\tau,e^{\sqrt{-1}(\theta+\vp+2k\pi)/m}}(\gamma)
-i_{\tau,e^{\sqrt{-1}(\theta+2k\pi)/m}}(\gamma)\right)\\
&=&\sum_{\omega^m=z}S^+_M(\omega).\ey
Similarly we have
$$ S^{-}_{M^m}(z)=\sum_{\omega^m=z}S^{-}_M(\omega). $$\qe

The following proposition gives the definition of the mean index.

\p (Theorem 1.5 of \cite{Lo4}) For any $\tau>0$ and
$\gamma\in\Pa_{\tau}(2n)$, there holds
\ee \hat i_{\tau}(\gamma)\equiv\lim_{k\to+\infty}
\frac{i_{k\tau}(\tilde\gamma)}{k}
=\frac{1}{2\pi}\int_{\Ua}i_{\tau,\omega}(\gamma)d\omega.\lb{756}\eee
In particular, $\hat i_{\tau}(\gamma)$ is always a finite real number,
which is called the {\bf mean index} per $\tau$ of $\gamma$.
\lb{758}\ep

\s{Variational properties of closed characteristics.}\lb{99chs3}

To solve the given energy problem (\ref{99che101}), we follow \S{V.3} of
\cite{Ek3}. Fix a given $\Sigma\in\Ha(2n)$ bounding a convex compact set
$C$. Let $j_C:\R^{2n}\to [0,+\infty)$ be the gauge function of $C$ defined
by
\ee j_C(0)=0\qquad {\rm and}\qquad j_C(x)=
\inf\{\lambda>0\,|\,\frac{x}{\lambda}\in C\}\quad
{\rm for}\quad x\ne 0. \lb{eq1.2}\eee
Fix a constant $\alpha$ satisfying $1<\alpha<2$ in this paper. As usual
we define the Hamiltonian function $H_{\Sigma,\alpha}:\R^{2n}\to [0,+\infty)$ by
\ee H_{\Sigma,\alpha}(x) = j_C(x)^{\alpha}, \qquad \forall x\in\R^{2n}.
\lb{eq1.3}\eee
Then $H_{\Sigma,\alpha}\in C^1(\R^{2n},\R)\cap C^2(\R^{2n}\bs\{0\},\R)$ is
convex and $\Sigma=H_{\Sigma,\alpha}^{-1}(1)$. It is well-known that
the problem (\ref{99che101}) is equivalent to the following problem
\ee \cases{\dot{z}(t) &$= JH_{\Sigma,\alpha}'(z(t)),  \qquad\forall t\in\R,$ \cr
z(1) &$=z(0). $\cr}\lb{99che301}\eee
Denote by $\J(\Sigma,\alpha)$ the set of all solutions $(\tau,x)$ of
the problem (\ref{99che301}) where $\tau$ is the minimal period of
$x$, and by $\tilde\J(\Sigma,\alpha)$ the set of all geometrically
distinct elements in $\J(\Sigma,\alpha)$. Note that elements in
$\J(\Sigma)$ and $\J(\Sigma,\alpha)$ are one to one correspondent
to each other.
The usual dual function $H_{\Sigma,\alpha}^{\ast}$ of $H_{\Sigma,\alpha}$ is defined by
$$ H_{\Sigma,\alpha}^{\ast}(x)=\sup_{y\in\R^{2n}}\{(x,y) -H_{\Sigma,\alpha}(y)\},   $$
where $(\cdot,\cdot)$ denotes the standard inner product of $\R^{2n}$.
For $1<\alpha<2$, let
\ee E_{\alpha}=\{u\in L^{(\alpha-1)/\alpha}
(\R/\Z),\R^{2n})\,|\,\int_0^1udt=0\}.\lb{99che302}\eee
The Clarke-Ekeland dual action functional
$f_{\alpha}:E_{\alpha}\to\R$ is defined by
\ee f_{\alpha}(u)=\int_0^1\{\frac{1}{2}(Ju,\Pi u)
+H_{\Sigma,\alpha}^{\ast}(-Ju)\}dt,\lb{99che303}\eee
where $\Pi u$ is defined by $\frac{d}{dt}\Pi u=u$ and $\int_0^1\Pi udt=0$.
Then $f_{\alpha}\in C^2(E_{\alpha},\R)$. Suppose $u\in E_{\alpha}\bs\{0\}$
is a critical point of $f_{\alpha}$. By \cite{Ek3}, there exists
$\xi_u\in\R^{2n}$ such that $z_u(t)=\Pi u(t)+\xi_u$ is a $1$-periodic
solution of the problem (\ref{99che301}).
Let $h=H_{\Sigma,\alpha}(z_u(t))$ and $1/m$ be the minimal period of $z_u$ for some
$m\in\N$. Define
\ee x_u(t) = h^{-1/\alpha}z_u(h^{(2-\alpha)/\alpha}t) \qquad
      {\rm and}\qquad \tau={1\over m}h^{(\alpha-2)/\alpha}.
\lb{99che304}\eee
Then there hold $x_u(t)\in\Sigma$ for all $t\in\R$ and
$(\tau,x_u)\in\J(\Sigma,\alpha)$. Note that the period $1$ of $z_u$
corresponds to the period $m\tau$ of the solution $(m\tau,x_u^m)$ of
(\ref{99che101}) with minimal period $\tau$.

On the other hand, every solution $(\tau,x)\in\J(\Sigma,\alpha)$
gives rise to a sequence $\{z^x_m\}_{m\in\N}$ of solutions of
the problem (\ref{99che301}), and a sequence $\{u^x_m\}_{m\in\N}$
of critical points of $f_{\alpha}$ defined by
\aa
  z^x_m(t) &=& (m\tau)^{-1/(2-\alpha)}x(m\tau t), \lb{99che307}\\
  u^x_m(t) &=& (m\tau)^{(\alpha-1)/(2-\alpha)}\dot{x}(m\tau t).\lb{99che308}
\eaa
For every $m\in\N$ there holds
\ee f_{\alpha}(u^x_m) = -(1-\frac{\alpha}{2})
      (\frac{2m}{\alpha}A(\tau,x))^{-\alpha/(2-\alpha)},
\lb{99che310}\eee
where
\ee A(\tau,x) = \frac{1}{2}\int_0^{\tau}(-J\dot{x}\cdot x)dt.
\lb{99che311}\eee

Following \S{V.3} of \cite{Ek3}, we denote by "$\ind$" the $S^1$-cohomology
index theory for $S^1$-invariant subsets of $E_{\alpha}$
defined in \cite{Ek3} (cf. also \cite{FR} of E. Fadell and P. Rabinowitz
for the original definition). For
$[f_{\alpha}]_c\equiv\{u\in E_{\alpha}\,|\,f_{\alpha}(u)\le c\}$
we define
\ee c_k = \inf\{c<0\,|\,\ind([f_{\alpha}]_c)\ge k\}. \lb{99che312}\eee
Then all these $c^k$'s are critical values of $f_{\alpha}$ and there hold
\aa
  -\infty &<& \min_{u\in E_{\alpha}}f_{\alpha}(u) = c_1 \le c_2\le \cdots
         \leq c_k \le c_{k+1}\le \cdots < 0, \lb{99che313}\\
  c_k &\to& 0 \qquad {\rm as}\quad k\to +\infty. \lb{99che314}\eaa

By Theorem V.3.4 of \cite{Ek3}, for each $k\in\N$ there exists a function
$u_k\in E_{\alpha}$ such that there hold
\aa
f_{\alpha}'(u_k) &=& 0 \qquad {\rm and} \qquad f_{\alpha}(u_k)=c_k,
\lb{99che315}\\
i_1^E(u_k) &\le& 2k-2 \le i_1^E(u_k)+\nu_1^E(u_k)-1. \lb{9che316}\eaa
Based upon this result and Proposition \ref{99chp203},
we have the following lemma.

\la \lb{99chl303} (cf. Lemma 3.1 of \cite{LZh2})
Suppose $^{\#}\tilde\J(\Sigma)<+\infty$, there exist an injection map
$p=p(\Sigma,\alpha)\colon\N\to\tilde\J(\Sigma,\alpha)\times\N$
such that for any $k\in\N$, $(\tau,x)\in\J(\Sigma,\alpha)$
and $m\in\N$ satisfying $p(k)=([(\tau,x)],m)$, there hold
\aa
f_{\alpha}^{\prime}(u^x_m) &=& 0 \quad {\rm and}
    \quad f_{\alpha}(u^x_m)=c_k, \lb{99che317}\\
i_{m\tau}(x^m) &\le & 2k-2+n\le i_{m\tau}(x^m)+\nu_{m\tau}(x^m)-1,
\lb{99che318}\eaa
where $u_m^x$ is defined by (\ref{99che308}).\el

\la\lb{99mll304} (cf. Corollary and 3.1 of \cite{LZh2})
Fix $\Sigma\in\Ha(2n)$ and $\alpha\in (1,2)$.
For any $(\tau,x)\in\J(\Sigma,\alpha)$ and $m\in\N$, there hold
\aa i_{(m+1)\tau}(x)-i_{m\tau}(x)&\ge&2,\lb{99che322}\\
i_{(m+1)\tau}(x)+\nu_{(m+1)\tau}(x)-1&\ge&i_{(m+1)\tau}(x)
>i_{m\tau}(x)+\nu_{m\tau}(x)-1,\lb{99che32}\\
\hat i_{\tau}(x)\ge 2.\lb{99che323}
\eaa
\el

\s{Proof of the main results.}\lb{99ml4}

Firstly we recall the common index jump theorem of \cite{LZh2}.

\t \lb{99cht401} (Theorem 4.3 of \cite{LZh2}) Let
$\gamma_k\in\Pa_{\tau_k}(2n)$, $k=1,\ldots,q$ be a finite collection
of symplectic paths. Let $M_k=\gamma(\tau_k)$ be the end points and
$\tilde\gamma_k$ be the iteration path of $\gamma_k$ defined by
(\ref{752}) for all $k=1,\ldots,q$. Denote by
$$ i_k^m=i_{m\tau_k}(\tilde\gamma_k)\quad\ox{and}\quad
\nu_k^m=\nu_{m\tau_k}(\tilde\gamma_k). $$
If $\hat i_{\tau_k}(x_k)>0$ for every $k=1,\ldots,q$, there exist
infinitely many $(N, m_1,\ldots,m_q)\in\N^{q+1}$ such that
\aa \nu_k^{2m_k-1}&=&\nu_k^1,\lb{99che423}\\
\nu_k^{2m_k+1}&=&\nu_k^1,\lb{99che424}\\
i_k^{2m_k-1}+\nu_k^{2m_k-1}&=&2N-(i_k^1+2S^+_{M_k}(1)-\nu_k^1),
\lb{99che425}\\
i_k^{2m_k+1}&=&2N+i_k^1,\lb{99che426}\\
i_k^{2m_k}&\ge&2N-n,\lb{99che427}\\
i_k^{2m_k}+\nu_k^{2m_k}&\le&2N+n\lb{99che428}\eaa
for all $k=1,\ldots,q$. Moreover, $M_k$ is elliptic if the equality in
one of the inequalities
of (\ref{99che427}) and (\ref{99che428}) holds.
\et

We need a lemma.

\la\lb{99mll401}Let $\gamma\in\Pa_{\tau}(2n)$ be a symplectic paths.
Let $M=\gamma(\tau)$ be the end point and $\tilde\gamma$ be the iteration
path of $\gamma$ defined by (\ref{752}). If $i_{\tau}(\gamma)\ge n$,
we have
\ee\lb{99mle501}i_{2\tau}(\tilde\gamma)
+2S^+_{M^2}(1)-\nu_{2\tau}(\tilde\gamma)\ge n.\eee\el

\bp For $\theta\in[0,2\pi)$, we define the splitting number
$S^-_{+,M}(e^{\sqrt{-1}\theta})$ by
\ee S^-_{+,M}(e^{\sqrt{-1}\theta})
=\nu_{\tau,e^{\sqrt{-1}\theta}}(\gamma)-S^-_M(e^{\sqrt{-1}\theta}).
\lb{99che421}\eee
By Corollary 4.13 of \cite{Lo4}, we obtain
\ee S^+_M(e^{\sqrt{-1}\theta})\ge 0. \lb{99che4191}\eee
By Lemma 6.3 of \cite{Lo4} and Lemma 4.2 of \cite{LZh}, we obtain
\aa \sum_{\theta\in (0,\pi)}S^-_M(e^{\sqrt{-1}\theta})
&+&S^-_{+,M}(1)+S^-_{+,M}(-1)\nn\\
&=& \sum_{\theta\in (0,\pi)}S^-_M(e^{\sqrt{-1}\theta})
+\frac{m(1)}{2}-r^-(1)+\frac{m(-1)}{2}-r^-(-1)\nn\\
&\le&\sum_{\theta\in (0,\pi)}q(e^{\sqrt{-1}\theta})
+\frac{m(1)}{2}+\frac{m(-1)}{2}\nn\\
&\le&n,
\lb{99che420}\eaa
where $q(e^{\sqrt{-1}\theta})$ is the negative Krein number of
$e^{\sqrt{-1}\theta}$ for $M$, $m(\pm 1)$ is the total multiplicity of
$\pm 1\in\sigma(M)$, $2r^-(\pm 1)$ is the number of Floquet
multipliers which arrive on the unit circle at ${\pm 1}$ for $M$.

Note that $S^+_M(\pm 1)=S^-_M(\pm 1)$ and $i_{\tau}(\gamma)\ge n$.
By Theorem 1.4 of \cite{Lo4} and Lemma \ref{99mll201}, we have
\y i_{2\tau}(\tilde\gamma)
+2S^+_{M^2}(1)-\nu_{2\tau}(\tilde\gamma)
&=&i_{\tau}(\gamma)+i_{\tau,-1}(\gamma)
+2S^+_M(1)+2S^+_M(-1)-\nu_{\tau}(\gamma)-\nu_{\tau,-1}(\gamma)\\
&=&2i_{\tau}(\gamma)+S^+_M(1)
+\sum_{\theta\in(0,\pi)}\left(S^+_M(e^{\sqrt{-1}\theta})
-S^-_M(e^{\sqrt{-1}\theta}\right)-S^-_M(-1)\\
& &+2S^+_M(1)+2S^+_M(-1)-\nu_{\tau}(\gamma)-\nu_{\tau,-1}(\gamma)\\
&=& 2i_{\tau}(\gamma)+2S^+_M(1)
+\sum_{\theta\in(0,\pi)}\left(S^+_M(e^{\sqrt{-1}\theta})
-S^-_M(e^{\sqrt{-1}\theta})\right)\\
& &-S^-_{+,M}(1)-S^-_{+,M}(-1)\\
&\ge& 2n-n\\
&=& n.\ey\qe

  Now we can give the proof of our main result. Fix $\Sigma\in\Ha(2n)$
and $\alpha\in(1,2)$. For any $(\tau,x)\in\J(\Sigma,\alpha)$, let
$\gamma_x\in\Pa_{\tau}(2n)$ be the fundamental solution of the linearized
system
$$ \dot y(t)=JH_{\Sigma,\alpha}''(x(t))y(t). $$
We call $\gamma_x$ the associated symplectic path of $(\tau,x)$.
Then the Maslov-type index $(i_{\tau}(x),\nu_{\tau}(x))$ and the
splitting number $S^+(x)$ of $x$ at $1$ are defined respectively by
\y i_{\tau}(x)&=&i_{\tau}(\gamma), \quad \nu_{\tau}(x)=\nu_{\tau}(\gamma),\\
S^+(x)&=&S^+_{\gamma(\tau)}(1).\ey

\la \lb{99mll402} Fix $\Sigma\in\SH(2n)$. For any
$(\tau,x)\in\J(\Sigma,\alpha)$, there holds $(\tau,-x)\in\J(\Sigma,\alpha)$,
and either $\O(x)=\O(-x)$ or $\O(x)\cap \O(-x)=\emptyset$. If
$\O(x)\cap \O(-x)\ne\emptyset$, we have
$$  x(t)=-x(t+\frac{\tau}{2}),  \qquad \forall t\in\R. $$
\el

\bp Since $\Sigma=-\Sigma$, there hold
\aa H_{\Sigma,\alpha}(x)&=&H_{\Sigma,\alpha}(-x),\lb{99mle430}\\
H^{'}_{\Sigma,\alpha}(x)&=&-H^{'}_{\Sigma,\alpha}(-x),\lb{99mle431}\\
H''_{\Sigma,\alpha}(x)&=&H''_{\Sigma,\alpha}(-x).\lb{99mle432}.\eaa
So for any $(\tau,x)\in\J(\Sigma,\alpha)$, there holds
$(\tau,-x)\in\J(\Sigma,\alpha)$.

If $\O(x)\cap \O(-x)\ne\emptyset$, there exist
$s_1,s_2\in[0,\tau]$ such that $x(s_1)=-x(s_2)$. By the fact
$x(\tau)=x(0)\ne 0$ and $x(t)\ne 0$ for any $t\in\R$, we have
$s_2-s_1\ne 0,\pm\tau$. Since $x(s_1+t)$ and $-x(s_2+t)$ satisfies the
same system
$$ \dot y=JH^{'}_{\Sigma,\alpha}(y), $$
we have $x(s_1+t)=-x(s_2+t)$ and hence $x(t)=x(2(s_2-s_1)+t)$ for $t\in\R$.
Since $\tau$ is the minimal period of $x$, we have $2(s_2-s_1)=\pm\tau$.
Therefore $x(t)=-x(t+\frac{\tau}{2})$. Specially there holds
$\O(x)=\O(-x)$. \qe

From $x(t)=-x(t+\frac{\tau}{2})$ we obtain $H''_{\Sigma,\alpha}(x(t))
=H''_{\Sigma,\alpha}(x(t+\frac{\tau}{2}))$.
Let $\gamma_x$ be the associated symplectic
path of $(\tau,x)$. Then we have
\ee \gamma_x(t+\frac{\tau}{2})=\gamma_x(t)\gamma_x(\frac{\tau}{2}),
\quad\forall t\in[0,\frac{\tau}{2}].\lb{99mle*}\eee

{\bf Proof of Theorem \ref{99mlt101}}\hspace{2mm}
It suffices to consider the case
$^{\#}\tilde\J(\Sigma)<+\infty$. By (\ref{99mle430}) and (\ref{99mle431}),
$(\tau,-x)\in\J(\Sigma,\alpha)$ if $(\tau,x)\in\J(\Sigma,\alpha)$.
Note that we have
$$(i_{m\tau_k}(x_k^m),\nu_{m\tau_k}(x_k^m))
=(i_{m\tau_k}(-x_k^m),\nu_{m\tau_k}(-x_k^m)),\qquad\forall m\in\N.$$
By Lemma \ref{99mll402}, we denote the elements
in $\tilde\J(\Sigma,\alpha)$ by
$$\tilde\J(\Sigma,\alpha)=\{[(\tau_j,x_j)]|j=1,\ldots,q_1\}
\cup\{[(\tau_k,x_k)],[(\tau_k,-x_k)]|k=q_1+1,\ldots,q_1+q_2\},$$
where $\O(x_j)=\O(-x_j)$ for $j=1,\ldots,q_1$, and
$\O(x_k)\cap \O(-x_k)=\emptyset$ for $k=q_1+1,\ldots,q_1+q_2$.
By Lemma \ref{99chl303}, we get an injection map
$p=p(\Sigma,\alpha)\colon\N\to\tilde\J(\Sigma,\alpha)\times\N$.
Note that $(\tau_k,x_k)$ and $(\tau_k,-x_k)$ have the same index intervals.
Thus by Lemma \ref{99chl303} we can further require that
$$ \im(p)\;\subset\{[(\tau_k,x_k)]\;|\;k=1,\ldots,q_1+q_2\}\times\N. $$

Set $i(k,m)=i_{m\tau_k}(x_k^m)$ and $\nu(k,m)=\nu_{m\tau_k}(x_k^m)$.
By Lemma \ref{99mll304} we have $\hat i_{\tau_k}(x_k)\ge 2$ for
$k=1,\ldots,q_1+q_2$.
Applying Theorem \ref{99cht401} to the associated symplectic path of
$$(\tau_1,x_1),\ldots,(\tau_{q_1+q_2},x_{q_1+q_2}),
(2\tau_{q_1+1},x_{q_1+1}^2),\ldots,(2\tau_{q_1+q_2},x_{q_1+q_2}^2),$$
we get infinitely many
$(N, m_1,\ldots,m_{q_1+2q_2})\in\N^{q_1+2q_2+1}$ such that
\aa \lb{99mle442}i(k,2m_k+1)&=&2N+i(k,1),\\
\lb{99mle443}i(k,2m_k-1)+\nu(k,2m_k-1)&=&2N-(i(k,1)+2S^+(x_k)-\nu(k,1)),\\
\lb{99mle444}i(k,2m_k)&\ge& 2N-n,\\
\lb{99mle445}i(k,2m_k)+\nu(k,2m_k)&\le& 2N+n\eaa
for $k=1,\ldots,q_1+q_2$ and
\aa\lb{99mle446} i(k-q_2,4m_k+2)&=&2N+i(k-q_2,2),\\
\lb{99mle447}i(k-q_2,4m_k-2)&+&\nu(k-q_2,4m_k-2)\nn\\
&=&2N-(i(k-q_2,2)+2S^+(x^2_{k-q_2})-\nu(k-q_2,2)),\\
\lb{99mle448}i(k-q_2,4m_k)&\ge& 2N-n,\\
\lb{99mle449}i(k-q_2,4m_k)+\nu(k-q_2,4m_k)&\le& 2N+n\eaa
for $k=q_1+q_2+1,\ldots,q_1+2q_2$.

  By (\ref{99mle444}), (\ref{99mle447}) and Lemma \ref{99mll401}, we have
\y i(k,2m_k)&\ge& 2N-n\\
&\ge&2N-(i(k,2)+2S^+(x_k^2)-\nu(k,2)),\\
&=& i(k,4m_{k+q_2}-2)+\nu(k,4m_{k+q_2}-2)\\
&>&i(k,4m_{k+q_2}-2).\ey
By (\ref{99mle445}), (\ref{99mle446}) and Lemma \ref{99mll401}, we have
\y i(k,2m_k)&<&i(k,2m_k)+\nu(k,2m_k)\\
&\le& 2N+n\\
&\le& 2N+i(k,2)\\
&=&i(k,4m_{k+q_2}+2)\ey
for $k=q_1+1,\ldots,q_1+q_2$. By Lemma \ref{99mll304} we have
$4m_{k+q_2}-2<2m_k<4m_{k+q_2}+2$. Thus
\ee m_k=2m_{k+q_2}, \qquad \forall k=q_1+1,\ldots,q_1+q_2. \lb{eq4.1}\eee

Denote by $([(\tau_{k(s)},x_{k(s)})],m(s))=p(N-s+1)$, where $s=1,\ldots,n$,
$k(s)\in\{1,\ldots,q_1+q_2\}$ and $m(s)\in\N$. By the definition of $p$
and (\ref{99che318}) we have
\aa i(k(s),m(s))
&\le& 2(N-s+1)-2+n=2N-2s+n   \nn\\
&\le& i(k(s),m(s))+\nu(k(s),m(s))-1.  \lb{eq4.2}
\eaa
So we have
\aa i(k(s),m(s))&\le& 2N-2s+n\nn\\
&<&2N+n\nn\\
&\le& 2N+i(k(s),1)\nn\\
&=&i(k(s),2m_{k(s)}+1), \lb{99mle451}
\eaa
for every $s=1,\ldots,n$, where we have used (\ref{99mle442}) in the
last equality. When $k(s)\le q_1$, by Lemma \ref{99mll401}, (\ref{99mle*}),
and (\ref{99mle443}), we have
\aa i(k(s),2m_{k(s)}-1)+\nu(k(s),2m_{k(s)}-1)
&=&2N-(i(k(s),1)+2S^+(x_{k(s)})-\nu(k(s),1))\nn\\
&\le&2N-n\nn\\
&\le& 2N-2s+n\nn\\
&\le& i(k(s),m(s))+\nu(k(s),m(s))-1. \lb{99mle452}\eaa
When $q_1<k(s)\le q_1+q_2$, by (\ref{eq4.1}) we have $m_{k(s)}=2m_{k(s)+q_2}$.
Then by (\ref{99mle447}), Lemma \ref{99mll401}, and (\ref{eq4.2}), we obtain
\aa i(k(s),2m_{k(s)}-2)+\nu(k(s),2m_{k(s)}-2)
&=&2N-(i(k(s),2)+2S^+(x^2_{k(s)})-\nu(k(s),2))\nn\\
&\le&2N-n\nn\\
&\le& 2N-2s+n\nn \\
&\le& i(k(s),m(s))+\nu(k(s),m(s))-1.\lb{99mle453}\eaa
By Lemma \ref{99mll304}, (\ref{99mle451}) and (\ref{99mle452})
we have $2m_{k(s)}-1<m(s)<2m_{k(s)}+1$ for $k(s)\le q_1$.
By Lemma \ref{99mll304}, (\ref{99mle451}) and (\ref{99mle453})
we have $2m_{k(s)}-2<m(s)<2m_{k(s)}+1$
for $q_1<k(s)\le q_1+q_2$.
Hence
\aa m(s)&=&2m_{k(s)},\qquad \ox{{\it if}}\; k(s)\le q_1 \lb{99mle454}\\
m(s)&\in&\{2m_{k(s)}-1,2m_{k(s)}\}\qquad \ox{{\it if}}\; q_1<k(s)\le q_1+q_2.
\lb{99mle455}\eaa
Since the map $p$ is injective, if there exist $s_1\ne s_2$ such that
$k(s_1)=k(s_2)\le q_1$, we must have $m(s_1)\ne m(s_2)$. Thus
$m_{k(s_1)}\ne m_{k(s_2)}$ by (\ref{99mle454}). This contradicts to
$k(s_1)=k(s_2)$. Similarly, if there exist $s_1\ne s_2$ such that
$k(s_1)=k(s_2)> q_1$, by the above arguments, we must have
$m(s_1)\ne m(s_2)$ and
$\{m(s_1), m(s_2)\}=\{2m_{k(s_1)}-1, 2m_{k(s_1)}\}$.
So there hold
$$ ^{\#}\{s\in\{1,\ldots,n\}|k(s)\le q_1\}\le q_1 $$
and
$$ ^{\#}\{s\in\{1,\ldots,n\}|k(s)>q_1\}\le 2q_2. $$
Therefore we have
\y ^{\#}\tilde\J(\Sigma)
&=& ^{\#}\tilde\J(\Sigma,\alpha)\\
&=&  q_1+2q_2\\
&\ge& n.\ey
This completes our proof. \qe

\bibliographystyle{abbrv}

\end{document}